\documentclass{article}
\usepackage[nonatbib,preprint]{nips_2018}

\usepackage[utf8]{inputenc} 
\usepackage[T1]{fontenc}    
\usepackage{hyperref}       
\usepackage{url}            
\usepackage{booktabs}       
\usepackage{amsfonts}       
\usepackage{nicefrac}       
\usepackage{microtype}      
\usepackage{graphicx}
\usepackage{amsmath}
\usepackage{verbatim}
\usepackage{enumitem}
\usepackage{amssymb}
\usepackage{amsthm}

\theoremstyle{theorem}
\newtheorem{theorem}{Theorem}[section]
\newtheorem{assumption}{Assumption}

\newcommand{\reals}{{\mbox{\bf R}}}

\newcommand{\diag}{\mathop{\bf diag}}

\newcommand{\dd}{\mathsf{d}}

\newcommand{\ie}{{\it i.e.}}

\title{On the Differentiability of the Solution to Convex Optimization Problems}

\author{
  Shane Barratt\\
  Department of Electrical Engineering\\
  Stanford University\\
  Stanford, CA 94305 \\
  \texttt{sbarratt@stanford.edu} \\
}

\begin{document}

\maketitle

\begin{abstract}
In this paper, we provide conditions under which one can take derivatives
of the solution to convex optimization problems with respect to problem data.
These conditions are (roughly) that Slater's condition holds,
the functions involved are twice differentiable, and that a
certain Jacobian matrix is non-singular.
The derivation involves applying the implicit function
theorem to the necessary and sufficient KKT system for optimality.
\end{abstract}

\section{Introduction}

Many engineering problems can be formulated as convex optimization problems and solved numerically.
Once an optimization problem is specified, it is common for each instance of the problem to be different, \ie, the objective function and constraints depend on situational problem data.
For example, in linear regression, the problem depends on the data matrix at hand.
It can be useful, in many cases, to characterize the
sensitivity of the solution (of the optimization problem being solved) to perturbations of the problem data.
The optimization problem might also depend upon
pre-specified hyperparameters, and in that case it can be useful to characterize the sensitivity of the solution to perturbations of the hyperparameters.

In this paper, we show how one can form the Jacobian matrix of the function that maps parameters
to a solution in a parametrized convex optimization problem.
This calculation requires several assumptions that are satisfied by many practical problems of interest.
We will also see that calculating the Jacobian can be fast compared to solving the actual optimization problem, provided one has already solved the problem.

\section{Preliminaries}

In this section, we review preliminaries: the setting of parametrized convex optimization,
the KKT conditions, and the implicit function theorem.
This section also serves as an introduction to notation used throughout the paper.
For a function $f:\reals^n \rightarrow \reals$,
we will denote its gradient by $\nabla_x f(x) \in \reals^n$
and for $f:\reals^n \rightarrow \reals^m$, we will denote its Jacobian by $\mathsf{D}_x f(x) \in \reals^{m \times n}$.

\subsection{Parametrized convex optimization}

We consider (convex) optimization problems of the form
\begin{equation}
\begin{aligned}
& \text{minimize}
& & f_0(x, \theta) \\
&&& f(x, \theta) \preceq 0 \\
&&& h(x,\theta)=0.
\end{aligned}
\label{eq:convexproblem}
\end{equation}
Here the vector $x \in \reals^n$ is the optimization variable of the problem,
the vector $\theta \in \reals^d$ is the problem data,
the function $f_0: \reals^n \times \reals^d \to \reals$ is the objective function,
the function $f: \reals^n \times \reals^d \to \reals^m$ forms the inequality constraints,
and the function $h: \reals^n \times \reals^d \to \reals^p$ forms the equality constraints.
The functions $f_i$ are assumed to be convex for fixed $\theta$, and the function $h$ is assumed to be affine for fixed $\theta$.
The optimal value $p^\star: \reals^d \rightarrow \reals$ of (\ref{eq:convexproblem}) for fixed $\theta$ is defined as
\begin{equation}
p^\star(\theta) = \inf \{f_0(x,\theta) | f(x,\theta)\preceq 0, h(x,\theta)=0\}.
\end{equation}
The set-valued \emph{solution mapping} $S(\theta): \reals^d \rightrightarrows \reals^n$ can be described as
\begin{equation}
S(\theta) = \{x \mid f(x,\theta)\preceq 0, h(x,\theta)=0, f_0(x,\theta)=p^\star(\theta)\}.
\end{equation}

We now state some assumptions that will be used in the sequel.
\begin{assumption}[Strong duality]
Slater's condition holds for (\ref{eq:convexproblem}).
\label{assumption:strongduality}
\end{assumption}
\begin{assumption}[Differentiability]
The functions $f_i$ are twice continuously differentiable in $x$, the function $f$ is continuously differentiable in $\theta$, the function $h$ is continuously differentiable in $\theta$ and $x$, and $\mathsf{D}_x f$ is continuously differentiable in $\theta$.
\label{assumption:differentiability}
\end{assumption}

\subsection{KKT conditions}

The following material is borrowed from Boyd and Vandenberghe~\cite{boyd2004convex}, with the only difference being that the objective and constraint functions are parametrized.
Define the Lagrangian
\begin{equation}
L(x,\lambda,\nu,\theta) = f_0(x,\theta) + \lambda^T f(x,\theta) + \nu^T h(x,\theta).
\end{equation}
Then necessary and sufficient optimality conditions for (\ref{eq:convexproblem}) are as follows.
The vector $\tilde{x}\in S(\theta)$ if and only if there are $(\tilde\lambda,\tilde\nu)$ that satisfy, with $\tilde{x}$, the KKT conditions:
\begin{equation}
\begin{aligned}
f(\tilde{x},\theta) &\preceq 0\\
h(\tilde{x},\theta) &= 0,\\
\tilde{\lambda}_i &\geq 0, \quad i=1,\ldots,m\\
\tilde{\lambda}_i f_i(\tilde{x},\theta) &= 0, \quad i=1,\ldots,m\\
\nabla_x L(\tilde{x},\tilde\lambda,\tilde\nu,\theta) &= 0.
\end{aligned}
\label{eq:kktconditions}
\end{equation}
Given a candidate solution $(\tilde{x},\tilde\lambda,\tilde\nu)$, let the set $G=\{i \mid \tilde\lambda_i=0 \; \text{and} \; f_i(\tilde{x}, \theta) = 0\}$.
We provide another assumption that is normally satisfied in practice, and allows us to ignore the first inequality in (\ref{eq:kktconditions}) when implicitly defining $\tilde{x}$.
\begin{assumption}[Emptiness of $G$]
The set $G=\emptyset$.
\label{assumption:emptiness}
\end{assumption}

\subsection{Implicit function theorem}

Consider functions $g: \reals^d \times \reals^n \rightarrow \reals^n$ and $x\in\reals^n$ implicitly defined by the equation $g(p,x)=0$ in which $p\in\reals^d$ acts as a parameter.
Define the \emph{solution mapping} $S:\reals^d \rightrightarrows \reals^n$ of this implicit equation as
\begin{equation}
S(p) = \{x \mid g(p,x)=0\}.
\label{eq:setmapping}
\end{equation}
Then the implicit function theorem is as follows, repeated from Dontchev and Rockafellar~\cite{dontchev2014implicit}.
\begin{theorem}[Implicit function theorem]
Let $g:\reals^d\times\reals^n\rightarrow\reals^n$ be continuously differentiable in a neighborhood of $(\bar{p},\bar{x})$ and such that $g(\bar{p},\bar{x})=0$, and let the partial Jacobian of $g$ with respect to $x$ at $(\bar{p},\bar{x})$, namely $\mathsf{D}_x g(\bar{p},\bar{x})$, be non-singular. Then the solution mapping $S$ defined in (\ref{eq:setmapping}) has a single-valued localization $s$ around $\bar{p}$ for $\bar{x}$ which is continuously differentiable in a neighborhood $Q$ of $\bar{p}$ with Jacobian satisfying
\begin{equation}
\mathsf{D}_p s(p) = -\mathsf{D}_x g(p, s(p))^{-1}\mathsf{D}_p g(p,s(p)) \;\; \text{for every} \;\; p \in Q.
\end{equation}
\end{theorem}
As is probably evident, we are going to apply the implicit function theorem to the KKT condition equations, which are a necessary and sufficient condition for optimality of (\ref{eq:convexproblem}).

\section{Finding the Jacobian}

We are going to reduce the KKT equations to an algebraic equation and apply the implicit function theorem.
We first let $z=(x,\lambda,\nu)$ for notational convenience and then define the function
\begin{equation}
g(z,\theta) = \begin{bmatrix}
\nabla_x L({x},\lambda,\nu,\theta) \\
\diag(\lambda) f({x},\theta) \\
h({x},\theta)
\end{bmatrix},
\end{equation}
where $\diag(\cdot)$ transforms a vector into a diagonal matrix.
If $g(\tilde{z},\theta)=0$ for some $\tilde{z}=(\tilde{x},\tilde\lambda,\tilde\nu)$ where $\tilde{x}$ and $\tilde{\lambda}$ are both feasible, and Assumptions~\ref{assumption:strongduality},~\ref{assumption:differentiability}, and~\ref{assumption:emptiness} are satisfied, then by (\ref{eq:kktconditions}) the vector $\tilde{x}$ is optimal.
Define the (partial) Jacobian
\begin{equation}
\mathsf{D}_{z} g(\tilde{z},\theta) = \begin{bmatrix}
\mathsf{D}_x \nabla_x L(\tilde{x},\tilde\lambda,\tilde\nu,\theta) & \mathsf{D}_x f(\tilde{x},\theta)^T & \mathsf{D}_x h(\tilde{x},\theta)^T \\
\diag(\tilde{\lambda}) \mathsf{D}_x f(\tilde{x},\theta) & \diag(f(\tilde{x},\theta)) & 0 \\
\mathsf{D}_x h(\tilde{x},\theta) & 0 & 0
\end{bmatrix}
\label{eq:partialx}
\end{equation}
and the (partial) Jacobian
\begin{equation}
\mathsf{D}_\theta g(\tilde{z},\theta) = \begin{bmatrix}
\mathsf{D}_\theta \nabla_x L(\tilde{x},\tilde\lambda,\tilde\nu,\theta) \\
\diag(\tilde\lambda) \mathsf{D}_\theta f(\tilde{x},\theta) \\
\mathsf{D}_\theta h(\tilde{x},\theta)
\end{bmatrix}.
\label{eq:partialalpha}
\end{equation}
Then, the following theorem holds.
\begin{theorem}[Differentiability of a Convex Optimization Problem]
If $g(\tilde{z},\theta)=0$, Assumptions~\ref{assumption:strongduality},~\ref{assumption:differentiability}, and~\ref{assumption:emptiness} hold, and $\mathsf{D}_x g(\tilde{z},\theta)$ is non-singular, then the solution mapping has a single-valued localization $s$ around $\tilde{x},\tilde\lambda,\tilde\nu$ that is continuously differentiable in a neighborhood $Q$ of $\theta$ with Jacobian satisfying
\begin{equation}
\mathsf{D}_\theta s(\theta) = -\mathsf{D}_z g(\tilde{x},\tilde\lambda,\tilde\nu,\theta)^{-1} \mathsf{D}_\theta g(\tilde{x},\tilde\lambda,\tilde\nu,\theta) \; \text{for every} \; \theta \in Q,
\label{eq:jacobian}
\end{equation}
with $\mathsf{D}_z g(\tilde{x},\tilde\lambda,\tilde\nu,\theta)$ defined in (\ref{eq:partialx}) and $\mathsf{D}_\theta g(\tilde{x},\tilde\lambda,\tilde\nu,\theta)$ defined in (\ref{eq:partialalpha}).
\end{theorem}

{\bf Proof.} The vectors $\tilde{x},\tilde\lambda,\tilde\nu$ are optimal if and only if $g(\tilde{z},\theta)=0$ and the formula for the Jacobian of the solution mapping (\ref{eq:jacobian}) follows from the implicit function theorem.

Interestingly, this theorem gets us the Jacobian of the optimal solution and dual variables.
Normally, however, we will just focus on the Jacobian with respect to the optimal solution.

Since most modern solvers perform primal-dual interior-point methods, they will at each step perform a factorization of the Hessian of the KKT matrix (\ref{eq:partialx}), which dominates the cost of calculating of the solution to (\ref{eq:jacobian}).
Thus, we can simply re-use this factorization to solve the system.

\section{Application to quadratic programs}
We now reconcile our results with those from Amos and Kolter~\cite{amos2017optnet}, where they derived the derivative of
the solution to Quadratic Programs (QPs) with respect to its input parameters, and used QPs as layers in a neural network.
They allude to the possibility of taking derivatives of general convex optimization problems, but focus on QPs.
They consider QPs of the form
\begin{equation}
\begin{aligned}
& \text{minimize}
& & \frac{1}{2}x^T Q(\theta) x + q(\theta)^T x \\
& \text {subject to}
&& G(\theta)x \preceq h(\theta) \\
&&& A(\theta)x=b(\theta).
\end{aligned}
\end{equation}
The partial Jacobian, using our framework, is just (suppressing the dependence on $\theta$ for convenience)
\begin{equation}
\mathsf{D}_x g(\tilde{x},\tilde\lambda,\tilde\nu,\theta)=
\begin{bmatrix}
Q & G^T & A^T \\
\diag(\tilde\lambda)G & \diag(G\tilde{x}-h) & 0\\
A & 0 & 0 \\
\end{bmatrix},
\end{equation}
which agrees with the left side of (6) in~\cite{amos2017optnet}.
Furthermore, the partial Jacobian for $\theta$ is
\begin{equation}
\mathsf{D}_\theta g(\tilde{x},\tilde\lambda,\tilde\nu,\theta) = \begin{bmatrix}
\dd Q\tilde{x}+\mathsf{D}_\theta q +\dd G^T\tilde \lambda+\dd A^T\tilde{\nu}\\
\diag(\tilde \lambda)(\dd G\tilde{x}-\mathsf{D}_\theta h)\\
\dd A\tilde{x}-\mathsf{D}_\theta b
\end{bmatrix}
\end{equation}
This matches the negative right-hand side of (6) in~\cite{amos2017optnet}.
Thus, our results confirm the findings of Amos and Kolter, generalize them to the general convex optimization setting, and provide conditions under which the Jacobian is well-defined.

\subsubsection*{Acknowledgments}
S. Barratt is supported by the National Science Foundation Graduate Research Fellowship
under Grant No. DGE-1656518

\small

\bibliographystyle{unsrt}
\bibliography{mybib}

\begin{thebibliography}{1}

\bibitem{boyd2004convex}
Stephen Boyd and Lieven Vandenberghe.
\newblock {\em {Convex Optimization}}.
\newblock Cambridge University Press, 2004.

\bibitem{dontchev2014implicit}
Asen~L Dontchev and R~Tyrrell Rockafellar.
\newblock {\em {Implicit Functions and Solution Mappings}}.
\newblock Springer-Verlag New York, 2014.

\bibitem{amos2017optnet}
Brandon Amos and Zico Kolter.
\newblock Optnet: Differentiable optimization as a layer in neural networks.
\newblock In {\em Proc. Intl. Conf. on Machine Learning (ICML)}, 2017.

\end{thebibliography}

\end{document}